\documentclass[12pt]{article}

\usepackage{amssymb,amsmath,amsfonts,amssymb}
\usepackage{graphics,graphicx,color,hhline}
\usepackage{bbm} 
\usepackage{euscript}
\usepackage{stmaryrd}
\usepackage[a-1b]{pdfx} 
\usepackage{mathtools}

\usepackage{authblk}

\def\@abssec#1{\vspace{.05in}\footnotesize \parindent .2in 
{\bf #1. }\ignorespaces} 

\graphicspath{{/EPSF/}{Figures/}}   

\newtheorem{theorem}{Theorem}[section]

\newtheorem{remark}[theorem]{Remark}

\def\ds{\displaystyle}
\def \Rm {\mathbb R}

\def \Cm {\mathbb C}

\def\un{{\mathbbmss{1}}}

\newcommand{\eps}{\varepsilon}

\newcommand{\be}{\begin{equation}}
\newcommand{\ee}{\end{equation}}
\newcommand{\bea}{\begin{eqnarray}}
\newcommand{\eea}{\end{eqnarray}}
\newcommand{\bee}{\begin{eqnarray*}}
\newcommand{\eee}{\end{eqnarray*}}

\newcommand{\bk}{\mathbf k}

\newcommand{\bx}{\mathbf x}

\newcommand{\bE}{\mathbf E}
\newcommand{\bB}{\mathbf B}

\def\fref#1{{\rm (\ref{#1})}}
\newcommand{\para}{\parallel}\newcommand{\bJ}{\mathbf J}
\newcommand{\boe}{\mathbf e}

\newcommand{\calQ}{\mathcal Q}

\newcommand{\calH}{\mathcal H}

\newcommand{\calF}{\mathcal F}

\newcommand{\calR}{\mathcal R}
\newcommand{\calG}{\mathcal G}

\newcommand{\calJ}{\mathcal J}

\newcommand{\calT}{\mathcal T}
\newcommand{\calK}{\mathcal K}

\newcommand{\cout}[1]{}

\begin{document}
\title{Time reversal of surface plasmons}
\author[,1]{Olivier Pinaud \footnote{pinaud@math.colostate.edu}}
 \affil[1]{Department of Mathematics, Colorado State University, Fort Collins CO, 80523}

 \maketitle

 \begin{abstract}
We study in this work the so-called ``instantaneous time mirrors'' in the context of surface plasmons. The latter are associated with high frequency waves at the surface of a conducting sheet. Instantaneous time mirrors were introduced in \cite{Fink-nphys}, with the idea that singular perturbations in the time variable in a wave-type equation create a time-reversed focusing wave. We consider the time-dependent three-dimensional Maxwell's equations, coupled to Drude's model for the description of the surface current. The time mirror is modeled by a sudden, strong, change in the Drude weight of the electrons on the sheet. Our goal is to characterize the time-reversed wave, in particular to quantify the quality of refocusing. We establish that the latter depends on the distance of the source to the sheet, and on some physical parameters such as the relaxation time of the electrons. We also show that, in addition to the plasmonic wave, the time mirror generates a free propagating wave that offers, contrary to the surface wave, some resolution in the direction orthogonal to the sheet. Blurring effects due to non-instantaneous mirrors are finally investigated.

   \end{abstract}

\section{Introduction}
This work is concerned with the concept of instantaneous time mirrors (ITM) recently introduced in \cite{Fink-nphys}. The main objective of this new technique is the control of waves by changing the underlying medium of propagation  as time evolves. In \cite{Fink-nphys}, waves at the surface of a water tank are perturbed by a sudden, strong shake of the tank. This results in the formation of a back-propagating wave (i.e. time-reversed) that spectacularly refocuses to re-create the initial source. At the mathematical level, ITM can be modeled  by singular in time perturbations of the constitutive parameters of the medium. A prototype for such models is the classical wave equation with velocity perturbed by a Dirac delta or one of its approximations. The problem was studied in \cite{BFP}, where the effect of the ITM is characterized in terms of an integral kernel that depends on the duration of the perturbation and on the modeling equation Green's function. Blurring is observed as the duration of the perturbation increases, and it is noteworthy to mention that the refocused wave is the time derivative of the initial source and not the source itself.  It is not completely direct to interpret the singular PDE, and an existence theory for the time singular wave equation was proposed in \cite{ITMPinaud}. Time-dependent media have been explored lately in different contexts, see \cite{lurie2007introduction,Milton-Mattei,ammari2021time,garnier2021wave,martin2020acoustics,borcea2019wave}.

One of the appeals of ITM is their relative experimental simplicity compared to classical time reversal. The latter involve the recording and reemission of the signal, see e.g. \cite{Fink-PT}, which could be difficult in practice.  ITM do not require such technical procedures, and open in principle the possibility to control quantum systems since ITM do not need measurements to generate the time-reversed wave \cite{QITM-D,QITM}. A central question at the core of ITM though is how to control the medium parameters. In \cite{Fink-nphys}, this is done by shaking the water tank, which changes the velocity of the surface water waves. Another experimental procedure is proposed in \cite{Low} and is based on the following observation: it is practically feasible to vary abruptly the density of available charge carriers in a graphene sheet, and this results in a singular time perturbation in the sheet conductivity. Such a perturbation acts as an ITM, and it becomes then possible to explore how the electromagnetic field generated by a dipole located close to the sheet is time-reversed by the ITM.

Our main objective in the present work is to characterize, in such a practical configuration, the point spread function (PSF) which describes the response of the perturbed system to an initial point source excitation. Another way to state the problem is  to ask whether it is possible to image the point source by using an ITM, and with which resolution. Wave propagation is modeled by the three dimensional time-dependent Maxwell's equations, and the conducting sheet (supposed to be flat) is taken into account via a jump condition on the horizontal (i.e. on the sheet's plane) magnetic field. The surface current on the sheet is obtained by Drude's model, which is an accurate description when doping in the sheet is sufficiently high. This system for wave propagation is standard, see e.g. \cite{MargLus1,low2017,Low,santosa,primerplasmon}. Note that while the equations are not time reversible due to the complex-valued conductivity of the sheet, the effects of irreversibility are mostly seen in a loss of amplitude of the time-reversed signal, offering therefore the possibility for sharp refocusing. Let us point out that the reference \cite{santosa} addresses a problem similar to the one considered here. Our work is different and complementary in that we fully characterize the time-reversed wave while \cite{santosa} establishes, among other facts, the generation of such wave in a simpler model without studying in detail the refocusing wave.

Physically, the emitted spherical wave interacts with the sheet, leading to the generation of a scattered wave that propagates freely after reflection and transmission, and a surface wave, referred to as the plasmonic wave, propagating along the sheet and evanescent in the orthogonal direction. An abrupt change in the sheet conductivity creates not only a back-propagating plasmonic wave, but also a back-propagating scattered wave. We will derive the PSF for these two waves. Due to its surface character, the plasmonic wave does not offer any out-of-plane resolution (i.e. range resolution), and we will see that the horizontal resolution (i.e. cross-range resolution) is proportional to the distance of the point source to the sheet, provided this distance is smaller than a characteristic length defined by some constitutive parameters of the sheet. One can then obtain excellent refocusing in the horizontal plane when the source is sufficiently close. We will also establish that the time-reversed scattered wave can provide some range resolution when the source is sufficiently far from the sheet. In such a case, the time-reversed plasmonic wave is negligible, and the scattered wave refocuses with a resolution that depends on the distance to the sheet and on some other parameters. We also investigate how the duration of the perturbation of the conductivity impacts refocusing, and show that some blurring is introduced the longer the perturbation.

The article is structured as follows: the model is introduced in Section \ref{setup}; ITM are defined in Section \ref{TR} for two types of perturbations, i.e. a Dirac delta and an approximation of it. Our main results are given in Section \ref{resu}: we derive expressions for the perturbed and unperturbed electromagnetic fields, which allow us to extract the time-reversed waves. We then obtain the PSF for the plasmonic and the scattered waves, and discuss their properties. Numerical simulations supporting our analysis are proposed. The details of the calculations are given in Section \ref{proofs}, and the article is ended with concluding remarks.


\paragraph{Acknowledgment.} This work is supported by NSF grant DMS-2006416.

\section{Setup} \label{setup}
We start with Maxwell's equations for the electromagnetic field $(\bE,\bB)$ in $\Rm^3$: with $\bx=(x,y,z)$, we have
\be \label{maxw}
 \frac{\partial \bB}{\partial t} +\nabla \times  \bE =0,\qquad \frac{1}{c^2} \frac{\partial \bE}{\partial t}- \nabla \times \bB + \mu_0 \bJ=0, \qquad \textrm{for } z \neq 0,
\ee
equipped with the initial conditions 
$$
\bB(t=0,\bx)=\bE(t=0,\bx)=0.
$$
Above, $c>0$ is the (constant) electromagnetic background velocity  and $\mu_0$ the permeability of free space. We suppose for simplicity that the background is the same below and above the conducting sheet located at $z=0$. Accounting for two different constant backgrounds  is possible at the price of more technicalities, but would not change qualitatively our results.

The current $\bJ$ is generated by a pointlike dipole located at $\bx_0=(0,0,z_0)$ and oriented along the vertical axis $\boe_z=(0,0,1)$, namely
\be \label{defJ}
\bJ(t,\bx)=U \delta(z-z_0)f(r) \delta(t) \boe_z, \qquad z_0>0, \quad r=\sqrt{x^2+y^2},
\ee
for $\delta$ the Dirac measure. With a ``genuine'' point dipole, we would have $f(r)=\delta(x)\delta(y)$. As discussed further in Section \ref{resu}, it turns out that with the typical experimental parameters, only sufficiently large horizontal (i.e. in the $(x,y)$ plane) wavenumbers in the source generate a plasmonic wave in the sheet. Smaller wavenumbers create surface evanescent waves located around the origin that do not propagate and decrease with time. There are two such modes, and while one has a large attenuation and is negligible, the other one has a weak attenuation and gives the leading contribution. When this latter mode is present, it is not possible to observe the time-reversed plasmonic wave associated with larger wavenumbers since its amplitude is weaker than that of the evanescent mode. We therefore suppose that the Fourier transform of $f$, with the convention
$$
\hat{f\;}(k)=\int_{\Rm^2} e^{-i (k_x x+k_y y)} f(r) dx dy, \qquad k=\sqrt{k_x^2+k_y^2},
$$
vanishes for $k\leq k_c$, where $k_c$ is a critical wavenumber for plasmonic waves to exist and which depends on the sheet's parameters. Its expression is given in Section \ref{resu}.

The orientation chosen for $\bJ$ simplifies the calculations as the system is invariant by rotation around the $z$ axis. For $\bB=(B_x,B_y,B_z)$, this leads to $B_z=0$ as is discussed in Remark \ref{sym}. The constant $U$ in \fref{defJ} is the amplitude of the emitted pulse, with physical unit Amp\`ere $\times $ meter $\times$ second (supposing $f$ is has the dimension of meter${}^{-2}$) to be consistant with \fref{defJ}.

We denote by $\llbracket \varphi \rrbracket(x,y)$ the jump of a function $\varphi$ across the plane $z=0$, that is
$$
\llbracket \varphi \rrbracket(x,y)=\lim_{z \to 0^+}\varphi(x,y,z)-\lim_{z \to 0^-}\varphi(x,y,z).
$$
The presence of the conductor sheet at $z=0$ is  modeled by the jump condition
\be \label{jump}
\llbracket \boe_z \times \bB\rrbracket=\mu_0 \bJ_s,
\ee
together with the continuity of $B_z$ at $z=0$. The surface current $\bJ_s \equiv \bJ_s(t,x,y)$ is found according to Drude's model:
\be \label{drude}
\frac{\partial \bJ_s}{\partial t}=-\frac{\bJ_s}{\tau}+D (t) \bE_\para(t).
\ee
Above, $\tau$ is the relaxation time of the electrons in the sheet, and $D(t)$ is the so-called Drude weight \cite{low2017}. The surface horizontal electric field $\bE_\para \equiv \bE_\para(t,x,y)$ in \fref{drude} is given by, with $\bE=(E_x,E_y,E_x)$,
$$
\bE_\para(t,x,y)= \big((\boe_z \times \bE) \times \boe_z \big)(t,x,y,0)=E_x(t,x,y,0) \boe_x+E_y(t,x,y,0) \boe_y.
$$
Note that the jump condition \fref{jump} can be included in \fref{maxw} by replacing $\bJ$ by $\bJ+\delta_{z=0} \bJ_s$. The sheet's conductivity relates $\bJ_s$ with $\bE_\para$, and depends therefore on $D(t)$. As explained in the Introduction, the time reversal of the plasmonic wave is accomplished by imposing an abrupt change in $D(t)$. Experimental values given in \cite{Low} for two different techniques show a switching time in the range $10^{-13}-10^{-14}$s, which is typically shorter than $\tau$ (of the order of $10^{-13}s$ \cite{primerplasmon}). For some time $T>0$, we will consider the case  $D(t)=D_0+ D_1(t-T)$, where $D_0$ is constant and
$$
D_1(t)=\alpha D_0 \delta(t), \qquad \textrm{or} \qquad D_1(t)=\alpha D_0 \chi_{\delta t}(t) \qquad \textrm{with} \qquad \chi_{\delta t}(t)=\frac{1}{\delta t} \chi\left(\frac{t}{\delta t}\right),
$$
where $\chi$ is a function with integral one supported in the ball of radius $1/2$  (we suppose as well without lack of generality that it is an even function to simplify some expressions further on). The function $\chi_{\delta t}$ is an approximation of $\delta(t)$ with $0<\delta t<T$ and $\alpha$ is a constant that has the dimension of time.

A mathematically rigorous existence theory of Maxwell's equations with the jump condition \fref{jump} and a singular $D(t)$ as above is beyond the scope of this work, and we will then only discuss informally why some objects are well-defined. We in particular implicitly assumed in \fref{jump}-\fref{drude} that $E_x$ and $E_y$ are continuous at $z=0$. This is easily seen as follows: take a smooth test function $\varphi\equiv\varphi(x,y)$ and integrate the left equation in \fref{maxw} over $\Rm^2 \times [-\eps,\eps]$. We find
\bee
\int_{- \eps}^{\eps} \int_{\Rm^2} \partial_t \bB\cdot \boe_x \, \varphi d \bx&=&-\int_{- \eps}^{\eps} \int_{\Rm^2} (\nabla \times \bE) \cdot \boe_x \varphi d \bx\\
&=&\int_{- \eps}^{\eps} \int_{\Rm^2} E_z \partial_y \varphi d \bx+\int_{\Rm^2} [E_y(x,y,\eps)-E_y(x,y,-\eps)] \varphi(x,y)dxdy.
\eee
Assuming $\partial_t B_x$ and $E_z$ are integrable, sending $\eps$ to zero shows that $E_y$ is continuous across the plane $z=0$. The continuity of $E_x$ is obtained similarly and by taking the scalar product with $\boe_y$ instead of $\boe_x$. 

\begin{remark} \label{sym}As mentioned earlier, we have $B_z=0$ in this configuration. This is explained below for completeness. For $z \neq 0$, we have indeed from \fref{maxw} that
\be \label{waveeq}
\partial_t^2 \bE-\Delta \bE=-\mu_0 \partial_t \bJ-\nabla (\nabla \cdot \bE).
\ee
The horizontal part to $\bJ$ is zero, and that of $\nabla (\nabla \cdot \bE)$ is a radial vector field since
$$
\nabla \cdot \partial_t \bE=-c^2 \mu_0 U f(r)\partial_z \delta(z-z_0), \qquad z \neq 0.
$$
Hence, \fref{waveeq} shows that the horizontal part to $\bE$ is also a radial vector field for $z \neq 0$. As a consequence, by \fref{maxw},
$$
\partial_t B_z=-\boe_z \cdot (\nabla \times \bE)=0, \qquad z \neq 0,
$$
leading to $B_z=0$ for $z \neq 0$ and actually everywhere by continuity of $B_z$ at $z=0$.
\end{remark}

The next section is dedicated to the ITM and time reversal. We decompose the field into various contributions and extract the ones giving the PSF.

\section{ITM and time reversal} \label{TR} We suppose as a start that $D_1(t)=\alpha D_0 \delta(t-T)$, namely that the ITM acts at $t=T$. Given Drude's equation \fref{drude}, such a choice for $D_1$ only makes sense when $\bE_\para$ is continuous with respect to the time variable. We will not prove this fact here and will only observe that since $\bE$ solves (in the distribution sense) for all $z$,
$$
\partial_t^2 \bE-\Delta \bE=-\mu_0 \partial_t \bJ-\mu_0 \delta_{z=0}\partial_t \bJ_s-\nabla (\nabla \cdot \bE),
$$
and since $\partial_t \bJ_s$ is as singular as  $\delta(t-T)$, the electric field $\bE$ has two more time derivatives than $\partial_t \bJ_s$, and is therefore continuous w.r.t. $t$. Hence, the surface current $\bJ_s$ is well-defined, and is discontinuous at $t=T$. One can find a proof of continuity in time for the solution to the  wave equation with delta-like coefficients in \cite{ITMPinaud}. It is quite technical, and its adaptation to our problem would require a separate work.

We  decompose $\bJ_s$ into
\bea \nonumber
\bJ_s(t)&=&\bJ_0(t) + \bJ_1(t):=D_0 \int_0^t e^{-(t-s)/\tau} \bE_\para(s)ds+  \int_0^t D_1(s) e^{-(t-s)/\tau} \bE_\para(s)ds\\
&=&\bJ_0(t)+\alpha D_0 \un_{t>T}\,e^{-(t-T)/\tau} \bE_\para(T). \label{expJs}
\eea
This shows that in order to obtain the perturbed solution for $t>T$, one simply needs to get the unperturbed solution (i.e. with $\alpha=0$) up to time $t=T$, and then solve the system with $\bJ_1$ as above and which is fully known. A time-reversed wave is created by the perturbation, and refocuses at $t=2T$. Since setting $\bJ$ as in \fref{defJ} is equivalent to setting $\bJ=0$ and
$$
\bE(t=0^+,\bx)=-c^2 \mu_0 U \delta(z-z_0) f(r)\boe_z,
$$
the PSF is obtained by considering $E_z$ at $t=2T$. We have then
$$
E_z(2T)=E_z^{(0)}(2T)+E_z^{(1)}(2T),
$$
where $E_z^{(0)}$ is the unperturbed solution. The perturbed solution $E_z^{(1)}$ is obtained by solving \fref{maxw} with $U=0$ and with $\bJ_s$ as in \fref{expJs} where $\bE_\para$ is replaced by $\bE_\para^{(0)}$, i.e. by the unperturbed surface horizontal electric field. We show in Section \ref{proofP} that the field $E_z^{(1)}(2T)$ can itself be decomposed into
\be \label{Ez1}
E_z^{(1)}(2T)=E_z^{(TR)}(2T)+E_z^{(F)}(2T)+E_z^{(M)}(2T),
\ee
where $E_z^{(TR)}$ corresponds to the time-reversed wave  and $E_z^{(F)}$ to the forward propagating wave. There are two contributions to $E_z^{(TR)}$: the purely plasmonic contribution, denoted by $E_z^{(P)}$, and the purely scattered wave one $E_z^{(S)}$. The term $E_z^{(M)}$ is a mixed plasmonic-scattered wave. It does not refocus since the dispersion relations of the plasmonic and the scattered wave are very different. The expressions of these terms above can be found in Section \ref{proofP}.

The values of $E_z(2T)$ around the emission point $(0,0,z_0)$ define the PSF. The contributions of the unperturbed solution $E_z^{(0)}$, of the forward propagating and mixed waves $E_z^{(F)}$ and $E_z^{(M)}$ are negligible compared to that of $E_z^{(TR)}(2T)$ since their dominating parts are supported away from $(0,0,z_0)$. We will therefore focus only on $E_z^{(TR)}$.

The situation is essentially the same in the regularized case where $D_1=\alpha D_0 \chi_{\delta t}$ with $\delta t \ll T$. We have now
$$
\bJ_1(t)=\int_0^t D_1(s-T) e^{-(t-s)/\tau} \bE_\para(s)ds,
$$
for $\bE_\para$ the total field, i.e. the sum of the unperturbed and perturbed fields.  Pick $s \in (T-\frac{\delta t}{2},T+\frac{\delta t}{2})$ (we recall that $\chi$ is supported in the ball of radius $1/2$). Then,
$$
\bE_\para(s)-\bE_\para^{(0)}(s)=\bE_\para(s)-\bE_\para(s-\delta t)+\bE_\para(s-\delta t)-\bE_\para^{(0)}(s).
$$
The second difference is  actually equal to $\bE_\para^{(0)}(s-\delta t)-\bE_\para^{(0)}(s)$ since $s-\delta t$ is not in the support of $\chi_{\delta t}$ and the perturbation has not occured yet. Hence, $\bE_\para(s)-\bE_\para^{(0)}(s)$ for $s$ as above is small in some sense around $s=T$  provided $\bE_\para$ and $\bE^{(0)}_\para$ have some regularity with respect to the time variable. It is shown in \cite{BFP,ITMPinaud}, in the context of the wave equation, that both fields have their time derivative bounded independently of $\delta t$. Again, adapting these proofs is beyond the scope of this work, and we will just claim the situation is similar here. We have indeed already observed at the beginning of the section that $\bE$ has two more time derivatives than $\partial_t \bJ_s$ (which is proportional to $\chi_{\delta t}$), and this implies that $\partial_t \bE$ behaves like the integral of $\chi_{\delta t}$ and is therefore bounded independently of $\delta t$. This yields $\bE_\para(s)-\bE_\para^{(0)}(s)=O(\delta t)$ for $s \in (T-\frac{\delta t}{2},T+\frac{\delta t}{2})$, and we then replace, as in the $D_1(t)=\alpha D_0\delta(t)$ case, $\bE_\para$ in $\bJ_1$ by the unperturbed field $\bE_\para^{(0)}$. The field $E_z^{(1)}$  verifies as a consequence the decomposition \fref{Ez1} with an additional error term of order $O(\delta t)$. We will see that the expressions of $E_z^{(P)}(2T)$ and $E_z^{(S)}(2T)$ are only slightly modified compared to the $D_1(t)=\alpha D_0 \delta(t)$ case.

\section{Results} \label{resu}

Our main results consist in the expressions of the time-reversed plasmonic and scattered waves $E_z^{(P)}(2T)$ and $E_z^{(S)}(2T)$, and in the analysis of their refocusing properties. Before stating those, we need to introduce some notation. Let first
\be \label{param}
\sigma_0=D_0 \tau, \qquad \eta=\frac{\mu_0 \sigma_0 c}{2}, \qquad \ell_0=\frac{\tau  c^2 \mu_0 \sigma_0}{2}=\eta c\tau.
\ee
The parameter $\eta$ is non-dimensional and typically small (we discuss experimental values further), $\ell_0$ has the dimension of a length and is referred to as the attenuation length, and $\sigma_0$ is the conductance of the sheet (in Siemens, or Ohms$^{-1}$). The quantity $\mu_0 c$ is the impedance of the surrounding medium measured in Ohms, and therefore indeed $\eta$ is non-dimensional. Let $\gamma=1-\eta^2$ and define $u_c>0$ by 
$$u^2_c=\frac{36\gamma-27-8\gamma^2+\sqrt{(36\gamma-27-8\gamma^2)^2+64\gamma^3(1-\gamma)}}{32}.$$
When $\eta \ll 1$, and therefore $\gamma \simeq 1$, we have $u_c\simeq 1/4$. For $k=\sqrt{k_x^2+k_y^2}$ defined earlier, consider the polynomial equation
\be \label{eqS}
s^4-2s^3+ \gamma s^2-k^2 \ell_0^2=0
\ee
that determines the plasmonic modes on the sheet, see Section \ref{proofunP}. When $k \ell_0 > k_c \ell_0:=u_c$, we show in Section \ref{proofunP} that the equation above has two complex conjugate roots that we denote by $s_\pm(k)=s_r(k)\pm i s_i(k)$, and two real roots that are not associated with physical solutions and are then ignored. When $k \ell_0 \leq  u_c$, there are only real roots and therefore no propagating plasmonic waves, and as explained in Section \ref{setup}, we only consider sources that have total horizontal momenta $k$ greater than $k_c$. Let moreover
$$
P'(s)=2-4 s-\frac{2 \eta^2}{1-s}.
$$
We need additional notations to define the time-reversed scattered wave. Let
\be \label{defT}
\sigma(\omega)=\frac{\sigma_0}{1+i \omega \tau}, \qquad \calT(\omega,k_z)=\frac{1}{1-\frac{\mu_0 \sigma(\omega) k_z c^2}{2 \omega}}.
\ee
Above, $\sigma(\omega)$ is the complex surface conductivity and $\calT(\omega,k_z)$ is the transmission coefficient of the sheet. For $|\bk|^2=k^2+k_z^2$ and $g$ a given function, we introduce
 \be \label{calH}
 \calH(T,k,z,g)=
 \int_\Rm e^{i cT |\bk| }e^{i k_z z } \calT(c|\bk|,k_z)g(c|\bk|)\frac{k_z d k_z}{|\bk|^2}.
 \ee
Let finally $$
  A=-\frac{2 \alpha \sigma_0(\mu_0 c^2)^2 U}{(4 \pi)^2}, \qquad
  \qquad \Delta_\para=\partial_{x^2}^2+\partial_{y^2}^2.
  $$

  We suppose that $cT>z_0$ for the emitted wave to have reached the sheet from the emission point. Our main theorem is the following:

  \begin{theorem} \label{th} Suppose $D_1(t)=\alpha D_0 \delta(t)$. Then, the time-reversed plasmonic and scattered waves admit the following expressions, for $z\geq 0$,
$$
E_z^{(P)}(2T,x,y,z)= A \Delta_\para \calJ_P(x,y,z), \qquad E_z^{(S)}(2T,x,y,z)= A \Delta_\para \calJ_S(x,y,z)
$$
where
  \bee
  \calJ_{P/S}(x,y,z)&=&\frac{1}{(2 \pi)^2}\int_{\Rm^2} e^{i (k_x x+k_y y)} \calK_{P/S}(k,z) \hat{f\,}(k) dk_x dk_y.
  \eee
 The kernels $\calK_P$ and $\calK_S$ are given by
  \bee
  \calK_P(k,z)&=& \frac{\Re\{(1-s_+)e^{i s_i (2 s_r-1) (z_0-z)/\ell_0}\}}{|P'(s_+)|^2} 16 \pi^2 e^{-2 s_r T/\tau}  e^{-(s_r-s_r^2+s_i^2) (z+z_0)/\ell_0}\\[3mm]
  \calK_S(k,z)&=& \Re \{\calH^*(T,k,z_0,1) \calH(T,k,z,\sigma/\sigma_0)\} .
  \eee
  Above, $\calH^*$ is the complex conjugate of $\calH$, $\Re$ denotes real part, and we recall that $s_+(k)=s_r(k)+is_i(k)$.
\end{theorem}

The proof of Theorem \ref{th} is given in Section \ref{proofs}. A first comment is, as in the case of the wave equation addressed in \cite{BFP}, that one does not directly recover a blurred version of the source, but here its horizontal Laplacian instead (up to multiplicative factors). Perfect refocusing corresponds to $\calJ_{P/S}(x,y,z)=\delta(z-z_0) f(r)$. With $\hat f(k)=\un_{k>\xi k_c}$ for some parameter $\xi>1$, we have $f(r)=\delta(x)\delta(y)-\int_{k\leq \xi k_c} e^{i (k_x x+k_y y)} dk_x dk_y/(2\pi)^2$, which behaves like $\delta(x)\delta(y)$ around $r=0$ since the second term is bounded. We study below the functionals $\calJ_P$ and $\calJ_S$ and quantify the amount of blurring compared to the perfect case $\calJ_{P/S}(x,y,z)\simeq \delta(z-z_0) \delta(x) \delta(y)$. The functional $\calJ_P$ is shown to peak in a region around $r=0$ that defines its horizontal resolution, while it decays exponentially away from the sheet and therefore does not localize around $z=z_0$. The functional $\calJ_S$ concentrates around the emission point $(0,0,z_0)$ in a region that defines its horizontal and vertical resolutions. We also discuss the regimes in which one functional dominates over the other.

\paragraph{Analysis of the plasmonic wave.} We suppose that $\hat{f \,}$ vanishes for $k\leq \xi k_c$ where $\xi$ is parameter such that $\xi>1$. When $\eta^2 \ll 1$, we show in Section \ref{proofunP} that $s_+$ is well approximated by
\be \label{approxS}
  s_+ \simeq \frac{1}{2}+\frac{\gamma-1}{4u}+ i \sqrt{u-u_c}, \qquad u=k \ell_0 > u_c.
  \ee
  The kernel $\calK_P$ then reads
  $$  \calK_P(k,z) \simeq \frac{\Re\{(1-s_+)e^{i (\gamma-1)\sqrt{u-u_c}  (z_0-z)/2u\ell_0}\}}{|P'(s_+)|^2} 16 \pi^2  e^{-T/\tau}  e^{-(z+z_0)/4\ell_0} e^{-(u-u_c)(z+z_0)/\ell_0}.
  $$
  The only term in $\calK_P$ that could potentially give some resolution in the vertical direction is the complex exponential in $z-z_0$. But on the one hand $\gamma-1 \simeq 0$, and on the other the phase vanishes for $u=u_c$ and in the limit of large $u$. The contribution of this term to $\calJ_P$ is therefore negligible, and $\calJ_P$ offers no vertical resolution, as expected. We then set the complex exponential to one in the sequel.

We turn now to the horizontal resolution. With \fref{approxS}, it follows that
  $$
   \frac{\Re\{(1-s_+)e^{i \sqrt{u-u_c} (\gamma-1) (z_0-z)/(2u\ell_0)}\}}{|P'(s_+)|^2}\simeq \frac{1}{32(u-u_c)}.
   $$
Set as an example $\hat f(k)=\un_{k>\xi k_c}$. After the change of variables $k \to (z+z_0)^{-1}k+\xi k_c$, we find
   \begin{align*}
     \calJ_P(x,y,z)\simeq \frac{\pi^2 e^{-\frac{T}{\tau}}e^{-\frac{(1+4(\xi-1)u_c)(z+z_0)}{4\ell_0}}}{2\ell_0 (z+z_0)} \int_0^\infty e^{- k} J\left(\left[k+\frac{\xi u_c(z+z_0)}{\ell_0}\right]\frac{r}{(z+z_0)}\right) \calQ(z,k) dk
   \end{align*}
   where
    \begin{align*}
    \calQ(z,k)= \frac{k+\frac{\xi u_c (z+z_0)}{\ell_0}}{k+\frac{u_c(z+z_0)}{\ell_0}[\xi-1]}.
   \end{align*}
  Above,  $J$ is $2 \pi$ times the zero-th order Bessel function of the first kind. The key quantity determining the horizontal resolution is $\zeta=\frac{\xi u_c(z+z_0)}{\ell_0}$. When $\zeta \ll 1$, then $\zeta$ can be ignored in the argument of $J$ in $\calJ_P$ to obtain simply $J(k r /(z+z_0))$, which goes to zero as $r \gg (z+z_0)$. Dominated convergence then shows that $\calJ_P$ tends to zero when $r\gg z+z_0$, showing that the horizontal resolution is $z+z_0$ when $\zeta \ll 1$. The term $\calJ_P$ can therefore be arbitrarily peaked around $r=0$ when $z=0$ and $z_0$ is close to the sheet. The term $\calQ$ has essentially no influence on the resolution since it tends to one for large $k$.

  When $\zeta \gg 1$, the dependency of $J$ in $k$ can be ignored and the Bessel function can be taken out of the integral. The term $\calJ_P$ is then proportional to $J(\xi u_c r /\ell_0)$, which shows that the horizontal resolution is now of order $\ell_0$ since $J(u) \to 0$ as $u \to +\infty$.

  For future comparisons with $\calJ_S$, we remark that $\calJ_P$ verifies when $\xi=2$,
  \be \label{caracP}
 \frac{\pi^3 e^{-T/\tau}e^{-(1+4(\xi-1)u_c)(z+z_0)/4\ell_0}}{\ell_0 (z+z_0)} \leq \calJ_P(0,0,z) \leq \frac{2\pi^3 e^{-T/\tau}e^{-(1+4(\xi-1)u_c)(z+z_0)/4\ell_0}}{\ell_0 (z+z_0)}.
 \ee

 \medskip
 
\begin{figure}[h!]
\centering
    \includegraphics[height=4.75cm, width=5.75cm]{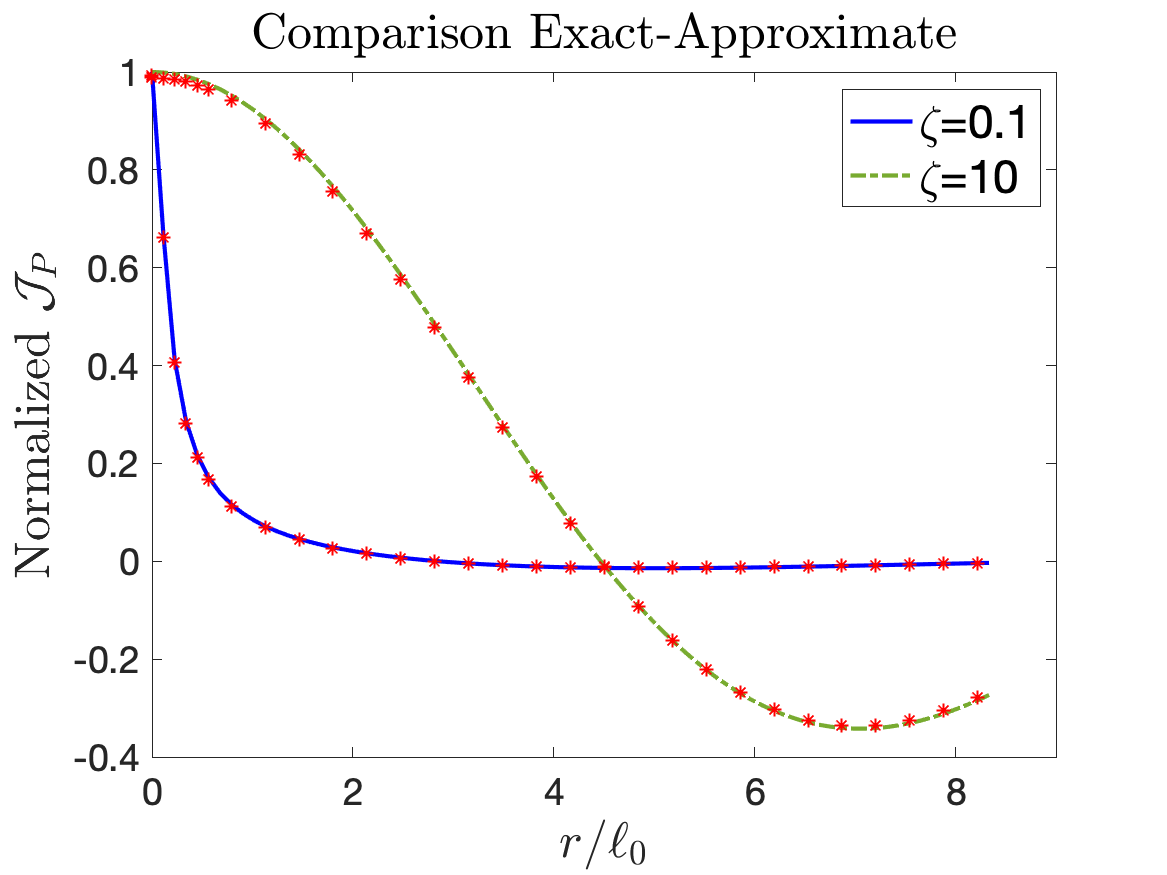} 
    \includegraphics[height=4.75cm, width=5.75cm]{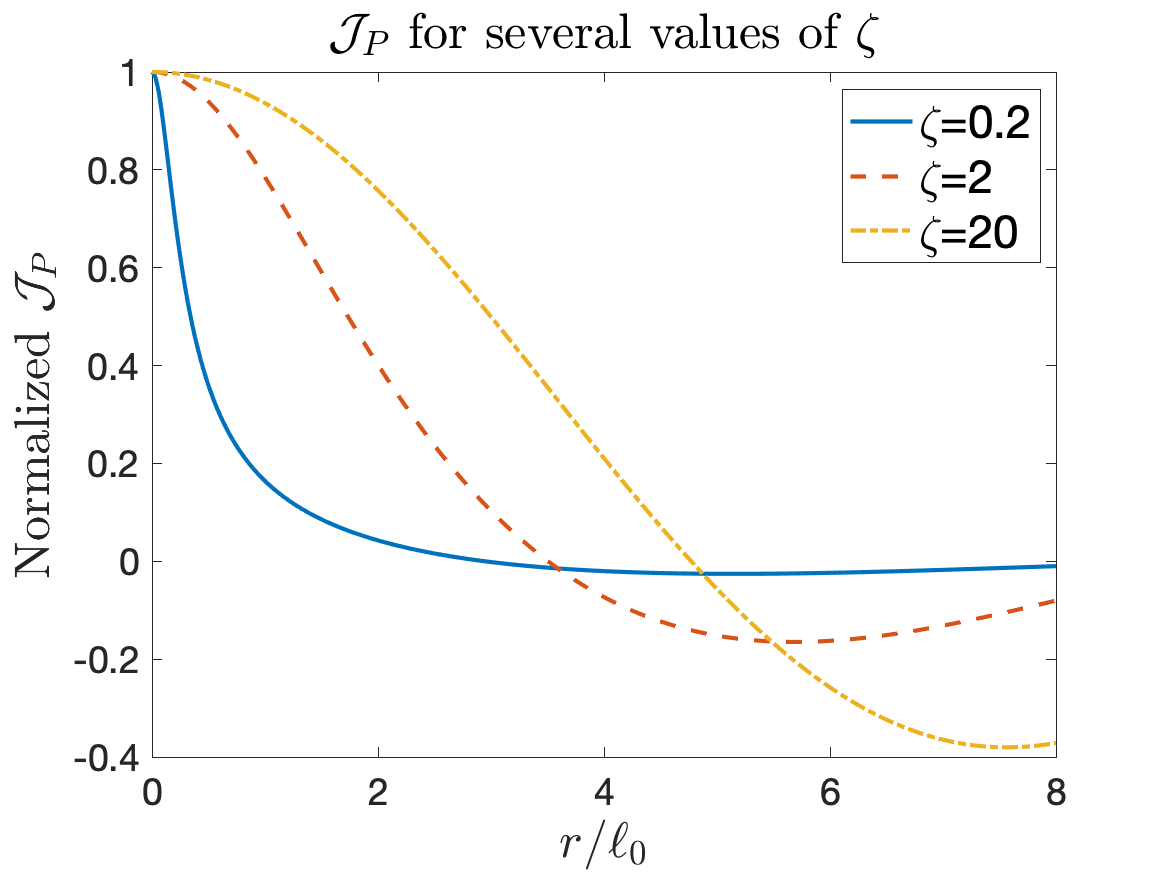}
  
\caption{Left panel: Comparison exact/asymptotic expressions of $\calJ_P$. The lines represent the exact expression, and the stars the approximate one. Right panel: representation of $\calJ_P$ for several $\zeta$. Observe the peaked behavior around $r=0$ for small $\zeta$.}
    \label{fig:JP}
  \end{figure}

 We compare in Figure \ref{fig:JP} the exact expression of $\calJ_S$ given in Theorem \ref{th} with the approximate one given above. Typical experimental values of the parameters are the following, see \cite{primerplasmon,low2017}: the drude weight $D$ is equal to $e^2E_F/\pi\hbar^2$, where $e$ is the electron charge, $\hbar$ the reduced Planck constant, and $E_F$ the Fermi energy. The latter is the tunable quantity as time varies, with maximal values of order of $0.4$eV. We then set e.g. $E_F=0.05$eV for the definition of $D_0$, with $E_F=0.4$eV corresponding to the value for the strong perturbation. With a  relaxation time $\tau$ of order $10^{-13}$s, this results in a graphene conductance $\sigma_0$ of the order of $6.10^{-4}$ Siemens, and assuming the surrounding medium has e.g. a refraction index of 2, we have $\mu_0 c \simeq 60 \pi $ Ohms, giving $\eta \simeq 0.1$ and $\gamma \simeq 0.99$. Moreover, since $\tau$ is of the order of $10^{-13}$ seconds, we find $\ell_0 \simeq 10^{-6}$ meters. For the calculation of the ``exact'' $\calJ_P$, we use numerical quadratures for the $k$ integral and find numerically the roots of \fref{eqS} for each $k$. We set $\xi=2$.

The exact and approximate functionals are represented in the left panel of Figure \ref{fig:JP} for $\zeta=0.1$ and $\zeta=10$, corresponding to the two asymptotic regimes described above. Observe the very good agreement. Both functionals are normalized using the value of the exact $\calJ_P$ at $r=0$. In the right panel of the figure, we represent $\calJ_P$ for $\zeta=0.2$, $\zeta=2$, and $\zeta=20$. For small $\zeta$, $\calJ_P$ is peaked around $r=0$ and offers therefore a very good resolution (of order $z+z_0$). For a larger $\zeta$, the resolution is limited to a few $\ell_0$, and we recognize the Bessel function in the case $\zeta=20$ as claimed earlier. 

  \medskip
  
  We now turn to the functional $\calJ_S$.
  
  \paragraph{Analysis of the scattered wave.}  We set as in the previous paragraph $\hat{f\,}(k)=\un_{k>\xi k_c}$ with $\xi=2$. The term $\calH$ is amenable to a stationary phase analysis under the assumption
  \be \label{asum}
  cT k \gg 1 \qquad \forall k \geq \xi k_c,
  \ee
 which is realized in most practical settings. Write indeed $cT= a z_0$, with $a>1$ for the emitted wave to reach the sheet. The smallest value of $k$ is $\xi u_c/ \ell_0$, and therefore $cT k \geq a\xi u_c z_0/\ell_0$. Since $\xi u_c \simeq 1/2$, the high frequency condition \fref{asum} becomes $a z_0 \gg 2 \ell_0$. This is immediately realized when the initial source if far away from the sheet, i.e. $z_0 \gg  \ell_0$, or when $T$ is sufficiently large in the case $z_0 \leq \ell_0$. Setting e.g. $z_0=\ell_0/5$ (which is $\zeta=0.2$ for $z=0$), gives $a \gg 10$, meaning $T$ must be greater than 10 times the time it takes for the initial pulse to reach the sheet.

  With
  $$
  \phi(z)=\left(1-\left(\frac{z}{cT}\right)^2 \right)^{1/2}, \qquad 
  \calT_0(k,z) = \frac{1}{1+\frac{\mu_0 \sigma_0 z}{2T[1+i c\tau k /\phi
      (z)]}},
   $$
   we show in Section \ref{proofP}, that under \fref{asum}, $\calH$ is approximated by, for $0<z<cT$, 
   \be \label{StatP}
   \calH(T,k,z,\sigma/\sigma_0) \simeq \left(
   \frac{2 \pi  z^2 }{\sigma_0^2 k (cT)^3 \phi(z) } \right)^{1/2} e^{i cT k\phi(z) } \calT_0(k,z)  \sigma(ck/\phi(z)).
 \ee
For $k \geq \xi k_c$, it turns out that $\calT_0(k,z) \simeq 1$ for $0<z<cT$, and $\calJ_S$ admits therefore the asymptotic expression:
 \begin{align} \label{AJS}
 &\calJ_S(z,x,y) \simeq  \sigma_0^{-1} C(z)
   \int_{\Rm^2} e^{i (k_x x+k_y y)} \Re \left\{e^{i cT k \Delta \phi(z)} \sigma(ck/\phi(z)) \right\} k^{-1}\hat{f\,}(k) dk_x dk_y,
 \end{align}
 where
 $$
 C(z)= \left(\frac{(2 \pi)^2  z^2 z_0^2 }{(cT)^6 \phi(z) \phi(z_0)} \right)^{1/2}, \qquad \Delta \phi(z)=\phi(z)-\phi(z_0).
 $$
 Remarking that

 $$
 \sigma(\omega)/\sigma_0=\frac{1}{1+(\tau\omega)^2}-\frac{i\tau \omega}{1+(\tau\omega)^2},
 $$
 $\calJ_S$ then reduces to
 $$
 \calJ_S(z,x,y) \simeq C(z)\int_{\xi k_c}^\infty J(kr) \left[ \frac{\cos[cT k \Delta \phi(z)]- \sin[cT k \Delta \phi(z)] c \tau k/ \phi(z)}{1+(c \tau k/\phi(z))^2}\right]dk,
 $$
 where we recall that $J$ is $2\pi$ times the zero-th order Bessel function of the first kind. When $z \neq z_0$, the fact that $\calJ_S$ is well-defined is directly established by integrating by parts the sine and cosine.

\begin{figure}[h!]
\centering
    \includegraphics[height=4.75cm, width=6.25cm]{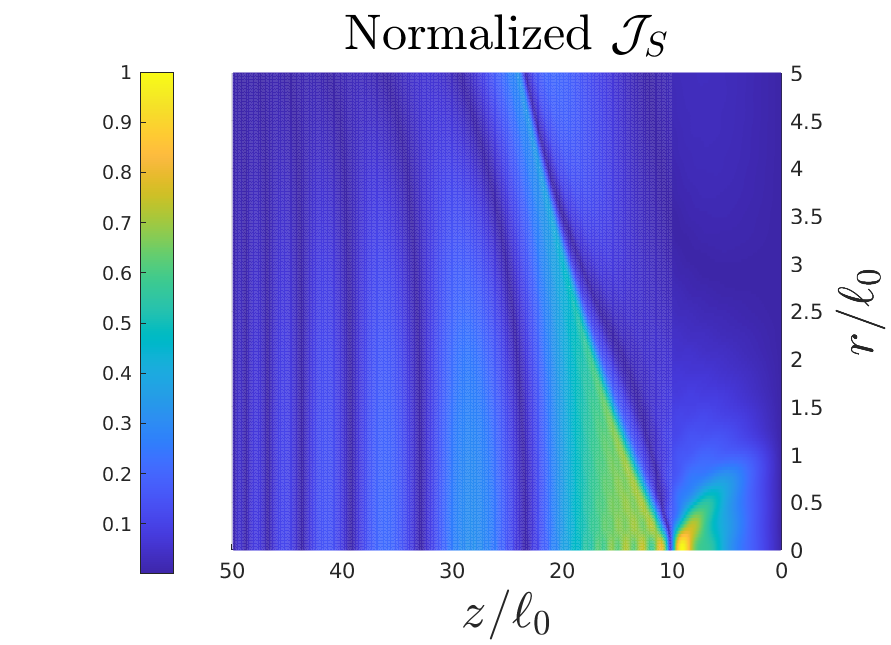} 
    \includegraphics[height=4.75cm, width=6.25cm]{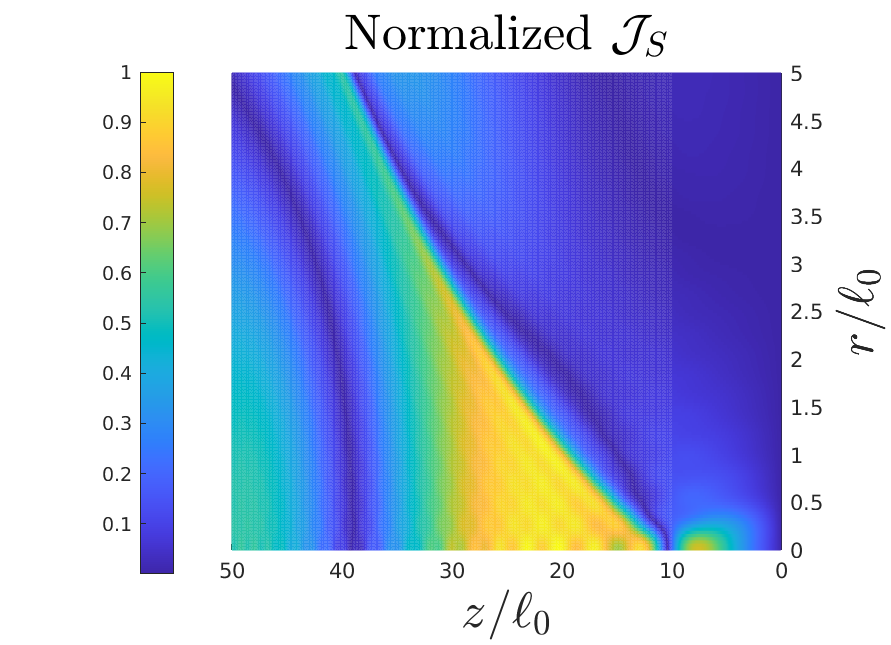}
  
\caption{Absolute value of $\calJ_P$ (normalized to one). Left panel: $cT=5z_0$, $z_0=10\ell_0$. Right panel: $cT=15z_0$, $z_0=10\ell_0$. Observe the loss of resolution with the increased $T$.}
    \label{fig:JS}
  \end{figure}

 We address now the refocusing properties of $\calJ_S$. Suppose first that $r \gg \ell_0 $. Then, $ kr \gg 1$ since $k \geq \xi u_c / \ell_0$, and
 \be \label{Basymp}J(kr) \simeq 2 \pi \sqrt{2/ \pi kr} \cos(kr-\pi/4).\ee
 Writing the integrand in $\calJ_S$ in terms of complex exponentials, it follows from the Riemann-Lebesgue Lemma that $\calJ_S\ll 1$  when
 $$
 |r\pm cT \Delta \phi(z)| \gg \ell_0.
 $$
 As a consequence, $\calJ_S$ is supported mostly in the region $|r\pm cT \Delta \phi(z)| \leq \ell_0$ when $r \gg \ell_0$. This is confirmed in Figure \ref{fig:JS}. Moreover, the asymptotic form of the Bessel function given above shows that $\calJ_S$ decreases as $1/\sqrt{r}$. This means that $\calJ_S$ is maximal in the region where $r \leq \ell_0$, resulting in a resolution of $\calJ_S$ in the horizontal plane of order $\ell_0$.

 Suppose now $r \ll \ell_0$. Then, $J(kr) \simeq 2 \pi$, and it is not difficult to see that $\calJ_S$ is small when 
 $$
  |cT \Delta \phi(z)| \gg \ell_0.
  $$
  When $(z/cT)^2 \ll 1$ and $(z_0/cT)^2 \ll 1$, this yields that $\calJ_S$ is maximal in the region where approximately
  \be \label{resoV}
  |(z+z_0) (z-z_0) | \leq 2 cT \ell_0.
  \ee
  When $z>z_0$, this condition is equivalent to $z-z_0\leq \sqrt{z_0^2+2cT\ell_0}-z_0$, showing that $\calJ_S$ peaks for $z> z_0$ in a region of vertical extent of order $\sqrt{z_0^2+2cT\ell_0}-z_0$. When $0\leq z \leq z_0$, \fref{resoV} is always verified when $2cT\ell_0 \geq z_0^2$. The latter condition is satisfied in the main regime of interest for $\calJ_S$   that we discuss in the next paragraph (where e.g. $cT=10 z_0$, $z_0=10 \ell_0$). Since $C(z) \to 0$ as $z \to 0$, it follows that $\calJ_S$ decreases away from $z=z_0$ as $z \to 0$. In summary, $\calJ_S$ is concentrated in a region around $(0,0,z_0)$ of horizontal extent of order $\ell_0$, of vertical extent of order $\sqrt{z_0^2+2cT\ell_0}-z_0$ above $z_0$, and decays to zero as $z\to 0$.

  In order to find an approximate peak value for $\calJ_S$, we have to be a little careful as setting naively $z=z_0$ yields the wrong result. We then first realize that the leading term in $\calJ_S$ when $r \ll \ell_0$ is the one proportional to the sine since $c \tau k$ is large. Then,
  $$
 \calJ_S(z,x,y) \simeq -\phi(z)C(z) /c\tau \int_{\xi k_c}^\infty \frac{\sin[cT k \Delta \phi(z)]}{k}dk.
 $$
 The last integral is written as
 $$
 \int_{\xi k_c}^\infty \frac{\sin[cT k \Delta \phi(z)]}{k}dk=\int_{0}^\infty \frac{\sin[cT k \Delta \phi(z)]}{k}dk-\int^{\xi k_c}_0 \frac{\sin[cT k \Delta \phi(z)]}{k}dk.
 $$
 The second term is negligible at $z \to z_0$, while the second one gives $\pm \pi /2$ depending on the sign of $\Delta \phi(z)$. It follows that a characteristic value for $\calJ_S$ in the region $r \ll \ell_0$ and $z \simeq z_0$ is
 \be \label{caracS}
 \calJ_S \simeq \frac{2\pi^3  z_0^2 }{(cT)^3 c\tau}.
 \ee
 Note that $\calJ_S$ changes sign around $z=z_0$, which is clearly observed in Figure \ref{fig:JS}.

We represent $|\calJ_S|$ (computed using \fref{AJS}) in Figure \ref{fig:JS} as a function of $(z,r)$. We only consider the case $z_0 \gg \ell_0$ since we will see below that this is the only relevant one for $\calJ_S$. We set e.g. $z_0=10 \ell_0$, with either $cT=5z_0$ (left panel) or $cT=15 z_0$ (right panel). We have as before $\eta=0.1$ and $\ell_0=10^{-6}$ meters. We observe as expected that the horizontal resolution is of the order of $\ell_0$, and that the vertical resolution above $z_0$ gets worse as $T$ increases. On the left panel, it is of the order of $\sqrt{z_0^2+2cT\ell_0}-z_0 \simeq  4\ell_0$ and the source is well-resolved, while $\sqrt{z_0^2+2cT\ell_0}-z_0=20 \ell_0$ on the right panel and there is a loss of resolution. 

  We now compare $\calJ_P$ and $\calJ_S$.

  \paragraph{Comparison $\calJ_P$-$\calJ_S$.} The best horizontal resolution for $\calJ_P$ is achieved at the sheet location $z=0$. When $z \simeq 0$, $\calJ_S \simeq 0$ and $\calJ_P$ dominates. We have already established that $\calJ_P$ does not offer any vertical resolution, and one can then ask if it is possible to obtain some with $\calJ_S$. When e.g. $z_0=\ell_0$ and $cT=15 z_0$, the ratio of the characteristic values of $\calJ_S$ and $\calJ_P$ at $z=z_0$ given in \fref{caracP}-\fref{caracS} is of the order of $10^{-3}$ when $r \leq \ell_0$. The ratio is even smaller when $z_0 \ll \ell_0$. This shows that $\calJ_P$ dominates and that the vertical location of the source cannot be determined when $z_0$ if of the order of $\ell_0$ or less. One could remark that the amplitude of $\calJ_P$ decreases at $T$ increases. Yet, for $\calJ_S$ to be larger than $\calJ_P$, one would need $T$ so large that $\calJ_S$ would not provide any vertical resolution since the latter is of the order of $\sqrt{z_0^2+2cT\ell_0}-z_0$. 

  The situation is different when $z_0 \gg \ell_0$, since it turns out that $\calJ_S$ dominates in that case in the vicinity of $z_0$. When e.g. $z_0=10 \ell_0$, the ratio is now of order $10^3$ for $cT=5 z_0$, and it becomes then possible to obtain some vertical resolution as described in the previous paragraph. In this regime, both $\calJ_S$ and $\calJ_P$ offer a horizontal resolution of order $\ell_0$

  We now consider the case when the perturbation $D_1$ is an approximation of a Dirac delta.

\paragraph{The regularized  case.} We assume again that $\hat f(k)=\un_{k>\xi k_c}$ with $\xi=2$. Calculations sketched in Section \ref{appdelta} show that the kernel $\calK_P$ becomes
$$
\calK_{P}(k,z) \hat \chi(2 \delta t\, s_i(k)/\tau),
$$
where we recall that $s_i$ is the imaginary part of the root $s_+$.
Above, $\hat \chi$ is the Fourier transform of $\chi$.  After the change of variables $k \to (z+z_0)^{-1}k+\xi k_c$ as before, we obtain the following for $\hat \chi$:
$$\hat \chi(2 \delta t\, \sqrt{k \ell_0/(z+z_0)+(\xi-1)u_c}/\tau).$$ Since, $\hat \chi(0)=1$, this shows that the regularization has essentially no effect on $\calJ_P$ when $\delta t \ll \tau \sqrt{(z+z_0)/\ell_0}$. In the opposite case when $\delta t \gg \tau \sqrt{(z+z_0)/\ell_0} $,
we find 
 \begin{align*}
     \calJ_P(x,y,z)\simeq \frac{\pi^2 \tau^2e^{-T/\tau}e^{-(1+4u_c)(z+z_0)/4\ell_0}}{2\ell_0^2 \delta t^2} \int_0^\infty J\left(\left[\frac{k(z+z_0) \tau^2 }{\ell_0 \delta t^2}+\frac{\xi u_c}{\ell_0}\right]r\right) \hat \chi (2 \sqrt{k})dk.
   \end{align*}
   The Bessel function is then approximately equal to $J(\xi u_c r/\ell_0)$ and can be taken out of the integral. This shows that the regularization introduced a loss of horizontal resolution that is now of order $\ell_0$ for all $z+z_0$, while it was $z+z_0$ before when $\zeta \ll 1$.

   Regarding $\calJ_S$, the kernel $\calK_S$ becomes
$$
\calK_{S}(k,z) \hat \chi(c \delta t\, k (\phi(z)+\phi(z_0))).
$$
When $(z/cT)^2 \ll 1$ and $(z_0/cT)^2 \ll 1$, we have
    $$
    \hat \chi(c \delta t\, k (\phi(z)+\phi(z_0)))\simeq \hat \chi(2 c \delta t\, k ).
    $$
    Since the minimal value of $k$ is $\xi u_c /\ell_0$, the regularization has no influence when $c \delta t \ll \ell_0$. In the intermediate case where $c \delta t =\ell_0$, the term $\hat \chi$ cannot be ignored but the regularization does not change qualitatively $\calJ_S$. When on the contrary $c \delta t \gg \ell_0$, we have $c \delta t\, k \gg 1$, leading to $\calJ_S \simeq 0$ since $\hat \chi(u) \to 0$ as $u \to +\infty$. The functional cannot then be used to obtain some vertical resolution in this case.
    




 \section{Proofs} \label{proofs}

We detail in this section the calculations leading to the results of Section \ref{resu}. We begin with some generalities, and then derive the perturbed and unperturbed solutions.
 
  \subsection{Generalities} We work in the Fourier space, and in order to compute Fourier transforms in time, we extend the fields by 0 for $t<0$. With the notation
  \be \label{nota}
  \hat{F}(\omega,k_x,k_y,z)=\int_{\Rm^3} e^{-i \omega t} e^{-i (k_x x+k_y y)} F(t,x,y,z) dtdxdy
  \ee
  for a given function $F$, Fourier transforming \fref{maxw} in all variables but $z$ yields
  \be \label{FourMax}
  \left\{
  \begin{array}{lll}
i \omega  \hat E_x &=& -c^2 \partial_z \hat B_y\\
i \omega \hat E_y &=& c^2 \partial_z \hat B_x\\
i \omega \hat E_z &=& i c^2 k_x \hat B_y-i c^2 k_y \hat B_x-c^2 \mu_0 U \hat{f\,}(k)\delta(z-z_0),
\end{array}\right.
                      \ee
                      as well as
      $$
  \left\{
  \begin{array}{lll}                
i \omega \hat B_x &=& \partial_z \hat E_y-i k_y \hat E_z\\
i \omega \hat B_y &=& i k_x \hat E_z -\partial_z \hat E_x\\
0 &=& i k_x \hat E_y - i k_y \hat E_x.
\end{array}\right.
$$
We recall that $B_z=0$ thanks to the rotational symmetry, which explains the last equation above. 
Easy algebra then gives the following equation for $\hat B_y$,
\be \label{eqBy}
\partial^2_z \hat B_y +\beta^2(\omega,k) \hat B_y = i k_x \mu_0 U \hat{f\,}(k) \delta(z-z_0), \qquad z \neq 0,
\ee
where
$$
\beta^2(\omega,k)=\frac{\omega^2}{c^2}-k^2, \qquad k^2=k_x^2+k_y^2.
$$
We will denote by $\beta(\omega,k)$ the complex square root of $\beta^2(\omega,k)$ with nonnegative imaginary part (and hence with branch cuts for $\Im(\omega)=0$, $(\Re(\omega))^2 \geq c^2k^2$). We now turn to the unperturbed solution.
\subsection{The unperturbed solution} \label{proofunP}

We set here $\alpha=0$ in the surface current $\bJ_s$ given in \fref{expJs} (i.e. $\bJ_1=0$). With the notation $\bJ_s=(J_{s,x},J_{s,y})$, \fref{eqBy} is equipped with the following jump condition at $z=0$
$$
\llbracket B_y\rrbracket=-\mu_0 J_{s,x},
$$
which follows from \fref{jump}, together with the continuity of $\partial_z \hat B_y$ at $z=0$ (which follows from the first equation in \fref{FourMax} and the continuity of $E_x$ at $z=0$ discussed in Section \ref{setup}). Solving the system for $\hat B_y$ gives
  $$
\hat B_y(\omega,k_x,k_y,z)=\frac{k_x \mu_0 U\hat{f\,}(k)}{2\beta(\omega,k)}\left\{
  \begin{array}{l}
   \ds \calR(\omega,k)e^{i \beta(\omega,k)(z+z_0)}+e^{i \beta(\omega,k)|z-z_0|}, \qquad z>0,\\
  \ds  \calT(\omega,k) e^{i \beta(\omega,k) (z_0-z)}, \qquad z<0
    \end{array}
  \right.
  $$
  with transmission and reflection coefficients $\calT$ and $\calR$ verifying
  $$
  \calT(\omega,k)=\frac{1}{1-\frac{\mu_0 \sigma(\omega) \beta(\omega,k)c^2}{2 \omega}}=1-\calR(k,\omega).
  $$
  Note that there is above a slight abuse of notation: to be consistent with the definition of the transmission coefficient $\calT$ given in \fref{defT}, we should write $\calT(\omega,\beta(\omega,k))$. We write here instead $\calT(\omega,k)$ for simplicity.
  With $\bE_\para^{(0)}=(E_x^{(0)},E_y^{(0)})$ the unperturbed surface horizontal electric field, we have then from \fref{FourMax} at $z=0$:
  $$
  \hat E^{(0)}_x(\omega,k_x,k_y,0)=-\frac{c^2}{i \omega} \partial_z \hat B_z=\frac{c^2 \mu_0 U \hat{f\,}(k)}{2\omega} e^{i \beta(\omega,k)z_0}\calT(\omega,k).
  $$
  After an inverse Fourier transform w.r.t. $\omega$, it follows that
  \be \label{expEx}
\calF E^{(0)}_x(T,k_x,k_y)=A_1\, k_x  \hat{f \,}(k)\int_{\Rm} e^{i \omega T} e^{i \beta(\omega,k) z_0 } \calT(\omega,k)  \omega^{-1} d\omega,
\ee
where $A_1=c^2 \mu_0 U /4 \pi$ and with the notation
$$
\calF E^{(0)}_x(t,k_x,k_y)=\int_{\Rm^2} e^{- i (k_x x+k_y y)}E^{(0)}_x(t,x,y) dx dy.
$$
The fact that the integral in \fref{expEx} is well-defined can be directly established using a stationary phase analysis since the function $\calT(\omega,k)  \omega^{-1}$ is locally integrable, smooth away from $\omega=ck$, and has no pole on the real axis as will be seen in the next paragraph.

We will only need $E_x^{(0)}$ in the sequel and not the $y$ component. The integral in \fref{expEx} can be decomposed into a branch contribution (due to the branch of the complex square root in $\beta(\omega,k)$), and a pole contribution due to the (complex) poles of $\calT$. The latter are solutions to
$$
1=\frac{\mu_0 \sigma(\omega) \beta(\omega,k)c^2}{2 \omega}.
$$
With $\omega= i s/\tau$ and the notations introduced in \fref{param}, this is equivalent to finding the solutions to \fref{eqS}. The latter are studied below.

  Note in passing that $B_y(t,x,y,z=0)=0$ when $z_0 > c t $ as the pulse has not reached the conducting sheet yet. This is seen for instance by observing that $\calR$ does not have poles or branch cuts at arbitrarily large imaginary values of $z$. Indeed, since $B_y$ is the space-time convolution of a function supported on the sphere of radius $ct-\sqrt{x^2+y^2+(z-z_0)^2}$ centered at $(0,0,z_0)$  and of the inverse Fourier transform of $\calR$, which is supported on $t \geq 0$, the result follows by simple inspection.

\paragraph{Analysis of the poles.} Let $u=k \ell_0$. We are interested in complex-valued solutions to \fref{eqS}, which are the ones associated with propagating modes. We will see that these exist only for $u$ sufficiently large when $\gamma \in [0,1]$ ($\gamma$ defined in \fref{param}). This is a simple consequence a standard formulas for quartic equations. Consider indeed the discriminant
$$
\Delta=16 u^2\big(-16u^4+(36\gamma-27-8\gamma^2)u^2+\gamma^3(1-\gamma)\big).
$$
When $\Delta<0$, \fref{eqS} has two distincts real roots and two complex conjugate roots. To check this condition, we consider the roots of $-16 u^2+(36\gamma-27-8\gamma^2)u+\gamma^3(1-\gamma)$. When $\gamma \leq 0$, both roots are negative and $\Delta$ is strictly negative when $u \neq 0$. When $\gamma \in (0,1]$, only one is positive, and $\Delta<0$ holds when $u^2$ is strictly greater than this positive root. A simple calculation shows that this root is equal to
$$u^2_c=\frac{36\gamma-27-8\gamma^2+\sqrt{(36\gamma-27-8\gamma^2)^2+64\gamma^3(1-\gamma)}}{32}.$$
When $u=u_c$, we have $\Delta=0$, and classical formulas for quartic equations show that there is a real double root (that turns into the complex roots when $u>u_c$) and two real simple roots.

Following this analysis, we will then only consider situations where $u > u_c$ for the propagatives modes to exist. We can actually say a little bit more about the roots. For $S_1$ and $P_1$ the sum and product of the complex roots, and $S_2$ and $P_2$ that of the other real roots, we can write \fref{eqS} as
$$
(s^2-S_1 s-P_1)(s^2-S_2 s-P_2)=x^4-(S_1+S_2)x^3+(S_1S_2-P_1-P_2)x^2+(S_1P_2+S_2P_1)x+P_1P_2.
$$
We find $P_1P_2=-u^2$ by identification with \fref{eqS}. Since $P_1>0$ when $u>u_c$, this shows that one of the real roots is positive and the other one is negative. Moreover, we have by inspection $S_1=-S_2 P_1/P_2$, and as a consequence $S_1$ and $S_2$ have the same sign. Since $S_1+S_2=2$, it follows that $S_1>0$, and that the complex roots have a positive real part. This shows that the propagative modes have a positive absorption coefficient and therefore decay exponentially.

\begin{figure}[h!]
\centering
    \includegraphics[height=3.75cm, width=3.75cm]{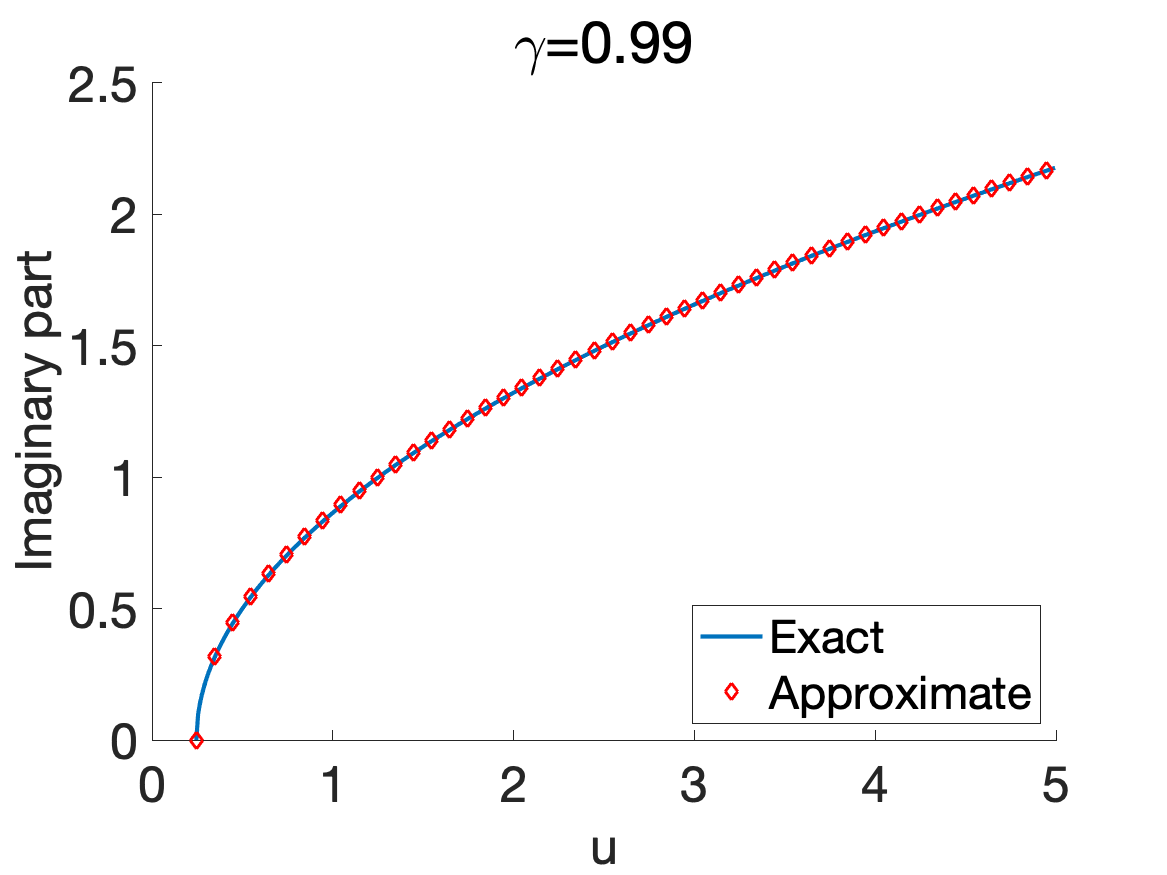} 
    \includegraphics[height=3.75cm, width=3.75cm]{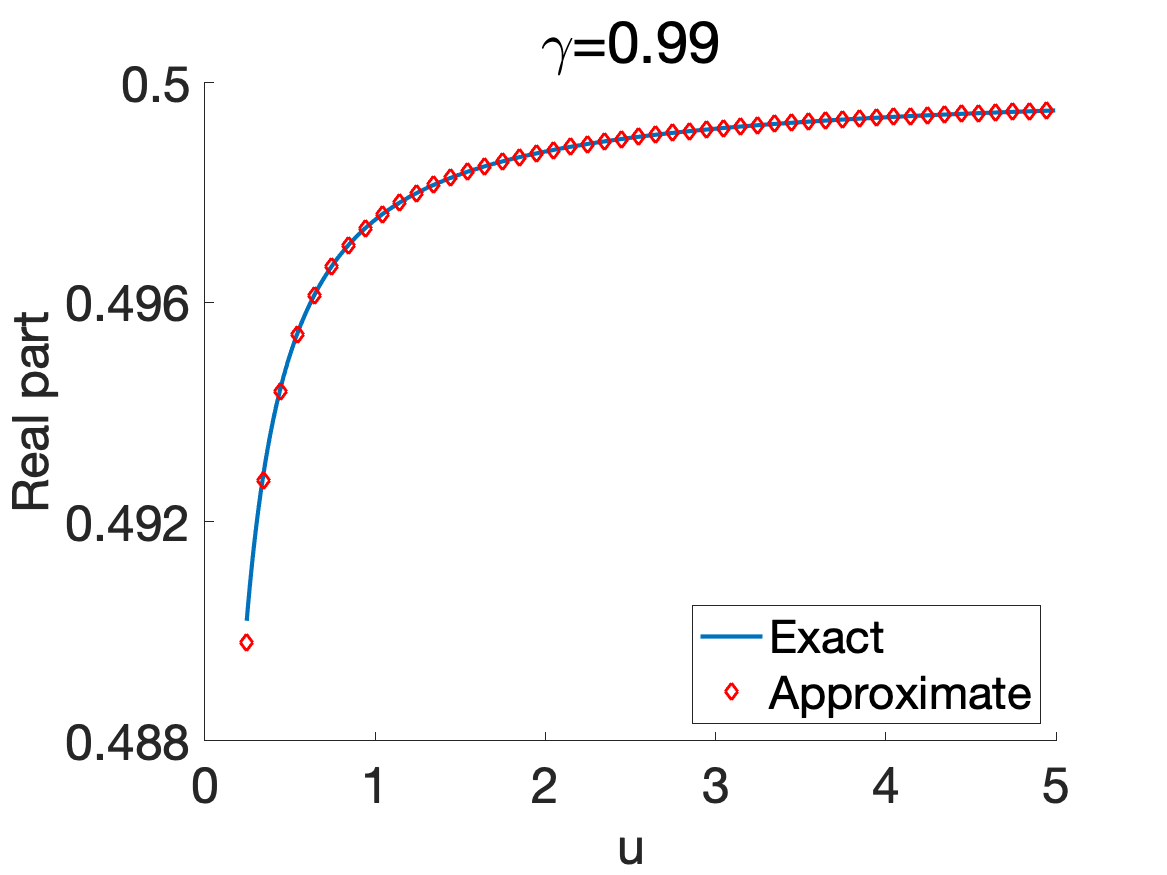}
     \includegraphics[height=3.75cm, width=3.75cm]{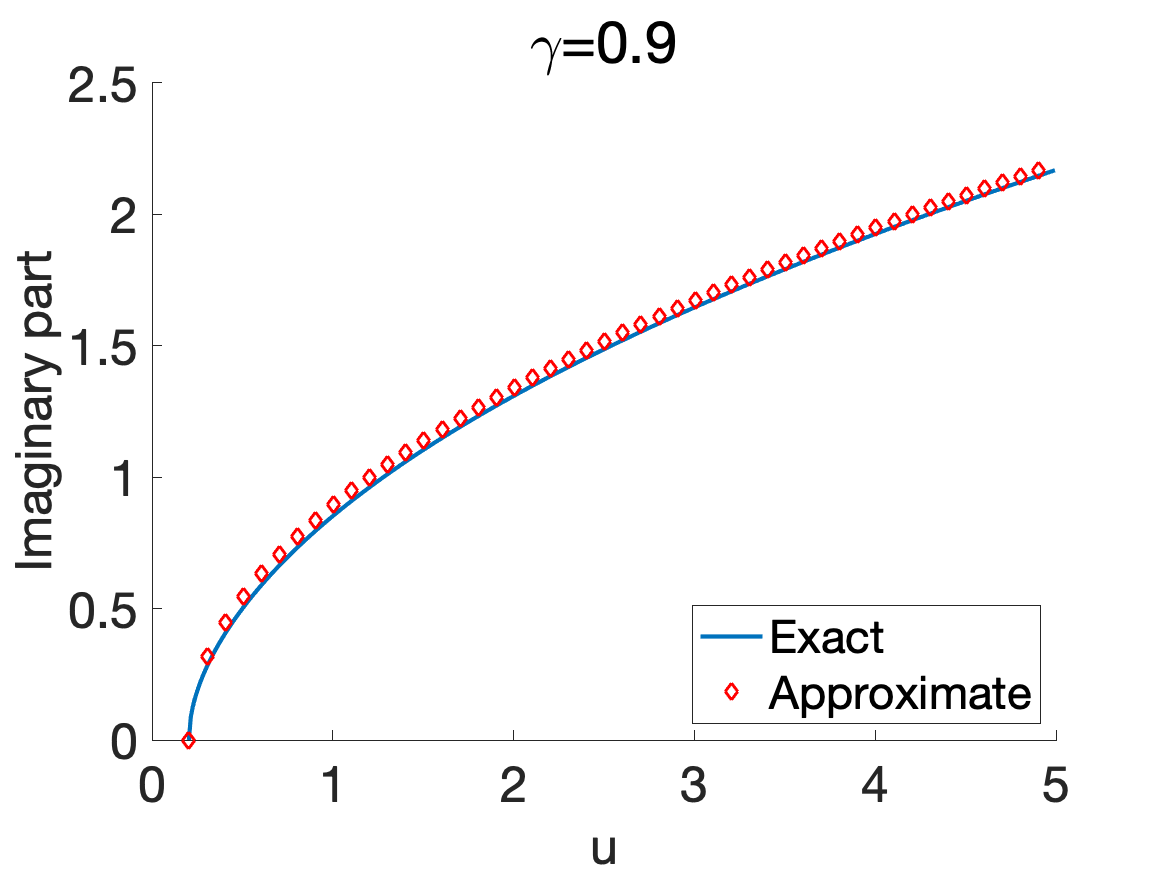} 
    \includegraphics[height=3.75cm, width=3.75cm]{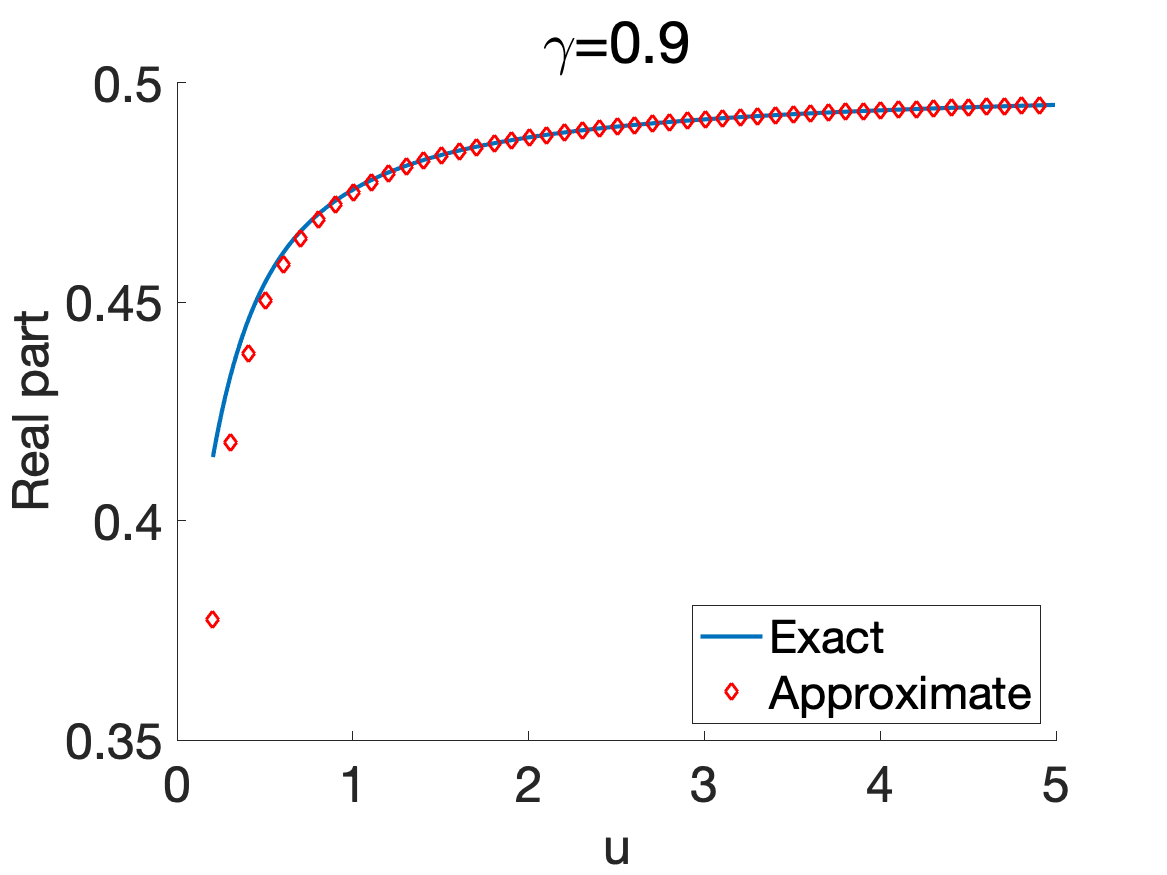}
\caption{Comparison exact/asymptotic expressions for the complex root $s_+$.}
    \label{fig:asymp}
\end{figure}

\textit{Asymptotics.} We investigate here the behavior of the complex roots for $u\gg1$. Write for this $s=a+ib$ for $a$ and $b$ real-valued. Separating real and imaginary parts in \fref{eqS} gives
\begin{align*}
  &b^3(2-4a)+b(4a^3-6a^2+2\gamma a)=0\\
  &b^4+(6a-6a^2-\gamma)b^2+a^4-2a^3+\gamma a^2-u^2=0.
  \end{align*}
  Assuming $u \gg 1$ and that $a$ remains bounded, a first crude estimation from the second equation above gives $b^4 \simeq u^2$, that is $b \simeq \sqrt{u}$. Expanding $a$ in power of $u^{-1}$ as $a=a_0+a_1u^{-1}+\cdots$, we find $a_0=1/2$, $a_1=(\gamma-1)/4$. We can then refine the estimate for $b$ as follows. Solving for $b^2$ in the second equation above gives
  $$
  b^2=\frac{1}{2}\left( \gamma+ 6 a^2 - 6 a +\sqrt{(\gamma+ 6 a^2 - 6 a)^2+4(u^2+2a^3-a^4-\gamma a^2)}\right).
  $$
  Setting $a=1/2$ and expanding in $u$ gives
  $$
  b^2 \simeq u+(\gamma-3/2)/2.
  $$
  Now, expanding $u_c$ around $\gamma=1$, we find $u_c=3/4-\gamma/2$, so that
  $$
  b \simeq \sqrt{u-u_c}.
  $$
  An approximation of the complex roots $s_{\pm}$ for $\gamma \simeq 1$ is then  given by
  $$
  s_\pm \simeq \frac{1}{2}+\frac{\gamma-1}{4u}\pm i \sqrt{u-u_c}, \qquad u > u_c.
  $$
  We compare the exact value and the approximation of the roots for $\gamma=0.9$ and $\gamma=0.99$ in Figure \ref{fig:asymp}. Note the excellent agreement even for small values of $u$ such that $u \simeq u_c$. 

We now transform \fref{expEx} using contour integration.

  \paragraph{Contour integrals.} First, we recall that $\beta(z,k)$ for $z \in \Cm$ has branch cuts on the real axis where $|z|\geq ck$. Consider then a function $F(z, \beta(z,k))$ holomorphic away from the branch cuts and from a finite number of poles with nonnegative imaginary parts. We suppose that $F$ converges sufficiently fast to $0$ on semi-circles of radius $R$ on the upper complex plane as $R \to \infty$ and consider
  $$
  I=\int_{\Rm} F(\omega, \beta(\omega,k)) d\omega.
  $$
  To avoid any confusion, it is best to see $\beta(\omega,k)$ as $c \beta(\omega,k)=\sqrt{\omega^2- i 0^+ \omega - c^2 k^2}$, where $i 0^+ \omega$ corresponds to an arbitrarily small absorption term added to \fref{maxw}. This shows that the integration in $I$ is done below the branch cuts on the real axis. Standard  contour integration then gives
  \begin{align} \nonumber
    I=I_{poles}&+\int_{c k}^\infty \big(F(\omega, -\beta(\omega,k))-F(\omega, \beta(\omega,k))\big) d\omega \\
    &+\int_{ck}^{\infty} \big(F(-\omega, \beta(\omega,k))-F(-\omega, -\beta(\omega,k))\big) d\omega, \label{contour}
  \end{align}
  where $I_{poles}$ is the pole contribution given by the residue theorem. We apply next these results to the calculation of $\calF E^{(0)}_x$.

\paragraph{The term $\calF E^{(0)}_x$.} With $F(\omega, \beta(\omega,k))=e^{i \omega T}e^{i \beta(\omega,k) z_0 } \calT(\omega,\beta(\omega,k)) \omega^{-1}$, the above analysis shows that $\calF E^{(0)}_x$ can be expressed as
$$
\calF E_x^{(0)}=\calF E_{x}^{(0,p)}+\calF E_{x}^{(0,s)},
$$
where $\calF E_{x}^{(0,p)}$ is the pole contribution (referred to as the plasmonic wave) and $\calF E_{x}^{(0,s)}$ the branch cut contribution (referred to as the scattered wave). For $z>0$, let
$$
\calG(t,z,s)= \frac{4 i \pi (1-s)}{P'(s)} e^{- s t/\tau}e^{-s(1-s) z/\ell_0}\qquad \textrm{with} \qquad P'(s)=2-4 s-\frac{2\eta^2}{1-s},
$$
for $s$ equal to the solutions to \fref{eqS} and $\eta$ defined in \fref{param}. Both real roots (positive and negative) lead to non viable physical solutions since they exhibit an exponential increase of $\calG$ (the positive pole becomes greater than $1$ as $k$ increases and therefore $s(1-s)<0$). They are then discarded. Since $P'(s_-)^*=P'(s_+)$, the residue theorem for $u>u_c$ yields the following expression for $\calF E_{x}^{(p)}$:

 \be \label{defEp}
 \calF E_{x}^{(0,p)}(T,k_x,k_y)= 2i A_1\, k_x\, \hat{f \;}(k) \Im\, \{\calG(T,z_0,s_+)\}.
 \ee
Regarding the scattered part, we first observe that
  $$
  \calT(\omega,\beta(\omega,k))=\calT(-\omega,-\beta(\omega,k))^*.
  $$
  After the change of variable $k_z=\beta(\omega,k)$ and some direct algebra based on \fref{contour}, we then find
 \be \label{defEs}
 \calF E_{x}^{(0,s)}(T,k_x,k_y)= -2i A_1\, k_x\, \hat{f \,}(k) \Im\, \{\calH(T,k,z_0,1)\},
 \ee
 where $\calH$ is defined in \fref{calH}.

 We turn now to the perturbed solution.
 
  \subsection{The perturbed solution} \label{proofP}
We set $U=0$, and solve Maxwell's equations with surface current $\bJ_s$ given by \fref{expJs}, and where $E_\para(T)$ is known and given by the $E_\para^{(0)}$ calculated in the previous section . With the notation \fref{nota}, we have
\be \label{J1}
\widehat{J}_{1,x}(\omega,k_x,k_y)=\alpha \sigma(\omega) e^{- i \omega T} \calF E^{(0)}_x(T,k_x,k_y).
\ee
 Solving Maxwell's equations gives the following expression for the perturbed magnetic field:
 $$
\hat B^{(1)}_y(\omega,k_x,k_y,z)=\frac{\mu_0}{2} \calT(\omega,k)\widehat{J}_{1,x}(\omega,k_x,k_y)\left\{
  \begin{array}{l}
   \ds  -e^{i \beta(\omega,k)z}, \qquad z>0,\\
  \ds   e^{-i \beta(\omega,k) z}, \qquad z<0
    \end{array}
  \right.
  $$
  Since \fref{FourMax} yields 
  $$
  \hat E_z^{(1)}=\frac{c^2 k^2}{\omega k_x} \hat B_y^{(1)},
  $$
  we find, for $z>0$, 
  $$
  \calF E^{(1)}_z(t,k_x,k_y,z)=-\frac{\alpha \mu_0 c^2 k^2}{4 \pi k_x} \calF E^{(0)}_x(T,k_x,k_y) \int_\Rm e^{i \omega (t-T)} e^{i \beta(\omega,k)z} \calT(\omega,k)    \sigma(\omega) \omega^{-1} d \omega.
  $$
  As in the previous section, the integral above can be decomposed into plasmonic and scattered parts. We find after direct calculations:
  \bee
  \int_\Rm e^{i \omega (t-T)} e^{i \beta(\omega,k)z} \calT(\omega,k)    \sigma(\omega) \omega^{-1} d \omega&=&2i  \Im\,\{ \sigma(i \tau^{-1} s_+) \calG(t-T,z,s_+)\}\\
  &&-2i\Im\,\{\calH(t-T,k,z,\sigma)\}.
  \eee
  Using the latter, together with \fref{defEp}-\fref{defEs}, we can then write at time $t=2T$, as announced in \fref{Ez1}, 
  $$
  \calF E_z^{(1)}=\calF E_z^{(TR)}+\calF E_z^{(F)}+\calF E_z^{(M)}=\calF E_z^{(P)}+\calF E_z^{(S)}+\calF E_z^{(F)}+\calF E_z^{(M)},
  $$
where $\calF E_z^{(TR)}$ is the time-reversed, backward propagating part of the perturbed wave, and $\calF E_z^{(F)}$ the forward propagating part. The former can be decomposed into a purely time-reversed plasmonic wave $\calF E_z^{(P)}$ and a purely time-reversed scattered wave $\calF E_z^{(S)}$. The term $\calF E_z^{(M)}$ is a mixed plasmonic-scattered wave. Their respective expressions are given by
  \begin{align*}
    &  \calF E_z^{(P)} =A_0 k^2\hat{f \,}(k)\, \Re \left\{\sigma(i \tau^{-1} s_+) \calG^*(T,z_0,s_+) \calG(T,z,s_+) \right\}\\
    &  \calF E_z^{(S)} =A_0 k^2\hat{f \,}(k)\, \Re \left\{\calH^*(T,k,z_0,1) \calH(T,k,z,\sigma)\right\},
  \end{align*}
and for the forward and mixed waves
  \begin{align*}
    & \calF E_z^{(F)} =-A_0 k^2\hat{f \,}(k)\, \Re \left\{ \big(\calG^*(T,z_0,s_+) +\calH(T,k,z_0,1)\big)\big( \sigma(i \tau^{-1} s_+)^* \calG^*(T,z,s_+)+\calH(T,k,z,\sigma\big))\right\}\\
      &  \calF E_z^{(M)} =A_0 k^2\hat{f \,}(k)\, \Re \left\{ \calG(T,z_0,s_+) \calH(T,z,s_+,\sigma)+ \sigma(i \tau^{-1} s_+)\calG(T,z,s_+)\calH(T,k,z_0,1)\right\}
  \end{align*}
  where
  $$
  A_0=\frac{2 \alpha (\mu_0 c^2)^2 U}{(4 \pi)^2}.
  $$
  The different definitions are motivated by the tendency of the dominating part of the wave to move from or to the point $r=0$: in $\calF E_z^{(P)}$ and  $\calF E_z^{(S)}$, the temporal phases compensate and yield backward propagating waves for $t \in [T,2T]$; in $\calF E_z^{(F)}$, the temporal phases add up (see the next paragraph for an approximate expression of $\calH$) to create a forward propagating wave with leading part supported away from $r=0$. This is seen by performing a standard stationary phase analysis for the inverse Fourier transform of $\calF E_z^{(F)}$. The term $\calF E_z^{(M)}$ is slightly different in that the temporal phases of $\calG$ and $\calH$ have different signs, but they do not compensate each other since one is significantly larger than then other. Indeed, using the asymptotic expressions for $s_+$ given in \fref{approxS} and that of the Bessel function given in \fref{Basymp}, the temporal phases in $\calF E_z^{(M)}$ read approximately, for $w=z,z_0$,
  $$
  i (cT \phi(w)\pm r) k-i \sqrt{k \ell_0-u_c}\, T/\tau, \qquad k \geq \xi \ell_0 u_c,
  $$
  which can be written as, 
  $$
  i cT/\ell_0 \left[(\phi(w)\pm r/(cT)) k \ell_0-\eta \sqrt{k \ell_0-u_c}  \right], \qquad k \geq \xi \ell_0 u_c.
  $$
  Since the term proportional to $\eta$ is small compared to the other one for large and small values of $k$, the mixed wave is localized around $r=cT \phi(w) \simeq cT$ when $(w/cT)^2 \ll 1$. This can be understood as follows. The mixed wave has two contributions: one originating from the plasmonic wave, and one from the scattered wave. The former creates a scattered wave at the time of the perturbation, and since the scattered wave propagates overall much faster than the plasmonic wave, it appears as if the scattered wave is created around $r=0$ at $t=T$. At time $t=2T$, the mixed wave is then localized  essentially at the same location as the generated scattered wave, which is around $r=cT$. There is a similar analysis for the plasmonic wave generated by the scattered wave. As a consequence, the mixed wave is indeed localized around $r=cT$ and its contribution can be neglected around the emission point.

We then focus on the time-reversed part $\calF E_z^{(TR)}$ for the study of the PSF since the contribution of $\calF E_z^{(F)}+\calF E_z^{(M)}$ is negligible in the vicinity of the source. 
  
We next simplify the expression of $\calH$ using stationary phase.

  \paragraph{Stationary phase analysis.} Rescaling $k_z$ as $k_z \to k k_z $, we find
 $$
 \calH(T,k,z,g)=
 \int_\Rm e^{i cTk( \sqrt{1+k_z^2}+zk_z/cT)} \calT(c k\sqrt{1+k_z^2},k k_z)g(ck \sqrt{1+k_z^2})\frac{k_z d k_z}{1+k_z^2}.
 $$
 The term $cTk$ is such that $cTk \gg 1$ according to assumption \fref{asum}, and the derivative of $\sqrt{1+k_z^2}+z k_z/cT$ vanishes when $z>0$ at the point $ k_z^*=-((\frac{cT}{z})^2-1)^{-1}$, where we recall that $cT>z$. With $g=\sigma$ or $g=1$, the term $\calT(c k\sqrt{1+k_z^2},k k_z)g(ck \sqrt{1+k_z^2})$ is smooth with bounded derivatives w.r.t. $k_z$, and a standard stationary phase procedure then yields the expression \fref{StatP}.

 We conclude this section with the case where the Dirac delta is regularized.


   \subsection{The regularized case.} \label{appdelta}
We only sketch the derivations here as the calculations are very similar to those of Section \ref{proofP}. The first difference is in the definition of $\widehat{J}_{1,x}$, that instead of \fref{J1} reads now
   $$
\widehat{J}_{1,x}(\omega,k_x,k_y)=\alpha \sigma(\omega) e^{- i \omega T} \int_{\Rm} \calF E^{(0)}_x(T+\delta t v,k_x,k_y) e^{-i \omega \delta t v} \chi(v)dv.
$$
This leads to the following expression of $ \calF E^{(1)}_z$:
\begin{align*}
  &\calF E^{(1)}_z(t,k_x,k_y,z)\\
  &=-\frac{\alpha \mu_0 c^2 k^2}{4 \pi k_x}  \int_{\Rm^2} \calF E^{(0)}_x(T+\delta t v,k_x,k_y) e^{i \omega (t-T-\delta t v)} e^{i \beta(\omega,k)z} \calT(\omega,k)    \sigma(\omega) \omega^{-1} \chi(v) dv d \omega.
\end{align*}
Proceeding as in Section \ref{proofP}, and using the fact that $\hat \chi$ is real and even since so is $\chi$, we obtain the expressions given in Section \ref{resu}.

\section{Conclusion}

We have studied in this work the time reversal of a plasmonic wave at the surface of a conducting sheet. Solving Maxwell's equations, we established the expression of the associated point-spread-function. On the one hand, we showed that the latter does not offer the possibility to image the vertical position of the source, as can be expected, and on the other that the resolution at which the horizontal location of the source can be determined depends on the distance $z_0$ of the source to the sheet: when $z_0 \ll \ell_0$, where $\ell_0$ is the attenuation length of the sheet, the resolution is of order $z_0$, while when $z_0 \gg \ell_0$, it is of order $\ell_0$. We also investigated the effects of the duration of the instantanenous time mirror on the point-spread-function, and quantified the amount of blurring introduced when the perturbation is not a Dirac delta. In addition to the plasmonic wave, we studied the time-reversed scattered wave created by the mirror. When $z_0\gg \ell_0$, the latter dominates over the plasmonic refocused wave and offers some vertical resolution.

This work raises a few natural questions. At the mathematical level, the well-posedness of Maxwell's system coupled to Drude's equation with a Dirac-type Drude weight remains to be done. Deriving uniform estimates in $\delta t$ in the case of the regularized delta is also of interest, and would justify rigorously the heuristic arguments we gave in Section \ref{TR}. Another question relates to control theory: we have seen that a refocusing wave is created by an ITM; by using other types of perturbations, which wave patterns can be generated? These questions will be investigated in future works.

\bibliographystyle{plain}
\bibliography{../bibliography.bib}
\end{document}